\documentclass[12pt, a4paper]{article}
\usepackage[latin1]{inputenc}
\usepackage[english]{babel}
\usepackage{indentfirst}
\usepackage{amstext}
\usepackage{setspace}
\usepackage{amsfonts}
\usepackage{textcomp}
\usepackage{amssymb}
\usepackage{amscd}
\usepackage{epsf}
\usepackage{graphicx}
\usepackage{epsfig}
\usepackage[mathscr]{eucal}
\usepackage{xr}
\usepackage{verbatim}
\usepackage{color}
\usepackage{shadethm}

\RequirePackage[reqno]{amsmath}

\newcommand {\demo}{\hskip -0.6cm{\bf Proof:  }}

\newcommand{\cqd}{\hfill{$\blacksquare$}}

\newcommand{\A}{\mathrm A}
\newcommand{\X}{\mathrm X}
\newcommand {\N}{\mathbb{N}}

\newcommand {\G}{\mathrm{G}}

\newcommand{\CC}{\mathrm C}
\newcommand{\U}{\mathrm U}

\newcommand{\K}{\ensuremath {\mathbb{K}} }

\newcommand{\Fo}{\ensuremath{\mathcal{F}_0(\X)}}
\newcommand{\Fa}{\ensuremath{\mathcal{F}_0(A)}}

\newcommand{\Fr}{\ensuremath{\mathcal{F}_0(R)}}

\newtheorem{theorem}{Theorem}[section]

\newtheorem{corolario}[theorem]{Corollary}
\newtheorem{definition}[theorem]{Definition}
\newtheorem{proposition}[theorem]{Proposition}

\newtheorem{obs}[theorem]{Remark}

\begin{document}

\title{Partial crossed products as equivalence relation algebras}
\maketitle
\begin{center}
{\large Viviane M. Beuter and Daniel Gonçalves}\\
\end{center}  
\vspace{8mm}

\abstract 
For a free partial action of a group in a set we realize the associated partial skew group ring as an algebra of functions with finite support over an equivalence relation and we use this result to characterize the ideals in the partial skew group ring. This generalizes, to the purely algebraic setting, the known characterization of partial C*-crossed products as groupoid C*-algebras. For completeness we include a new proof of the C* result for free partial actions.
\doublespace

\section{Introduction}

The groupoid approach to C*-algebras given by Jean Renault in \cite{Renault} is one of the main concepts in the modern theory of operator algebras, as an ever growing number of C*-algebras may be realized and studied as groupoid C*-algebras (AF algebras, Cuntz and Cuntz-Krieger algebras, tilings C*-algebras are just a few we can mention, see \cite{Abadie, Gonc, Gonc1}). Forteen years after the work of Renault, the notion of partial actions was introduced by Exel (see \cite{Exel1}) and McClanahan (see \cite{Mc}), and, as many C*-algebras were expressed and studied as partial crossed products (AF algebras, Cuntz-Krieger algebras of infinite matrices and Bunce-Deddens algebras are a few examples we can mention, see \cite{Abadie}), the theory grew in importance. In 2004 Abadie (see \cite{Abadie}) established that every partial crossed product of a commutative C*-algebra can be seen as a groupoid C*-algebra. In 2005 Dokuchaev and Exel (see \cite{Ex}) started the study of purely algebraic partial actions and their associated partial skew group rings, providing a new insight into the theory of partial actions, which was followed by a broading of the knowledge in the field, see \cite{Avilaferrero, Ex, Ex3, GR1, Ferrerolazzarin, BR} for example.

It is interesting to note that many results in the theory of C* partial crossed products have equivalent versions in the purely algebraic setting and, as the theory develops, the interaction between the areas increase. Although the study of partial skew group rings is still underdeveloped when compared with its C*-counterpart, we can mention a few examples where the interaction between the areas seems to be benefitial to both. For example proposition 2.1 of \cite{ExGiGo}, where partial actions of countable groups over second countable, compact spaces whose envelope space is Hausdorff are characterized can be obtained, see remark 2.2 of \cite{ExGiGo}, from the algebraic version described in \cite{Ex}. Also, conditons for simplicity of skew group rings and applications to topological dynamics (and hence to the associated C*-algebras) have been studied in \cite{Oinert, Oinert1}. Some of these results have recently been generalized to partial skew group rings, with applications to partial actions on compact sets (see \cite{Gonc2, GR2}). 

In this paper we intend to give purely algebraic versions of some known results in the theory of C*- partial crossed products. In particular, to give the reader some motivation, we start with a simpler, and algebraic in flavour, proof of Abadie´s characterization of C*-partial crossed products as groupoid algebras, for the case of free partial actions of countable groups acting on unital commutative C*-algebras. In this case the groupoid can be seen as an étale equivalence relation and, as a consequence, the construction of the groupoid C*-algebra falls within reach of a much bigger audience. In the second part of the paper we generalize the well know relation between partial dynamical systems and C*-partial dynamical systems to the purely algebraic level. We then proceed to generalize, to the algebraic level, Abadie´s result, that is, we show that the
partial skew group ring associated to a free algebraic partial action on a set is isomorphic to an algebra of functions with finite support over an equivalence relation. We finish the paper showing how to use the characterization just mentioned to obtain a 1-1 correspondence between ideals and $R$-invariant subsets. In light of recent results characterizing Leavitt path algebras as partial skew group rings, see \cite{GR1}, it is interesting to note that one can use this last correspondece to derive the known ideal structure (see \cite{tomforde}) of Leavitt path algebras associated to finite graphs with no cycles. Before we proceed we recall for reader´s convenience some key definitions below.

\begin{definition} A partial action of a group $\G$ on a set $\Omega$ is a pair $\theta=(\{\Delta_{t}\}_{t\in \G},$ $\{h_{t}\}_{t\in \G})$, where for each $t\in \G$, $\Delta_{t}$ is a subset of $\Omega$ and $h_{t}:\Delta_{t^{-1}} \rightarrow \Delta_{t}$ is a bijection such that:
\begin{enumerate}
	\item $\Delta_{e} = \Omega$ and $h_{e}$ is the identity in $\Omega$;
	\item $h_{t}(\Delta_{t^{-1}} \cap \Delta_{s})=\Delta_{t} \cap \Delta_{ts}$;
	\item  $h_{t}(h_{s}(x))=h_{ts}(x),$ $x \in \Delta_{s^{-1}} \cap \Delta_{s^{-1} t^{-1}}.$
\end{enumerate}

If $\Omega$ is a topological space, we also require that each $\Delta_{t}$ is an open subset of $\Omega$ and that each $h_{t}$ is a homeomorphism of $\Delta_{t^{-1}}$ onto $\Delta_{t}$. 

Analogously, a pair $\theta = (\{ D_{t} \}_{t \in \G} , \{ h_{t} \}_{t \in \G} )$ is a partial action of $\G$ on an algebra $\A$ if each $D_{t}$ is a closed two sided ideal and each $h_{t}$ is an isomorphism of $D_{t^{-1}}$ onto $D_{t}$. In case $\A$ is a C*-algebra we also require that each  $h_{t}$ is a *-isomorphism.
\end{definition}

\begin{definition} A partial action $(\{\Delta_{t}\}_{t\in \G},$ $\{h_{t}\}_{t\in \G})$ is said free if, for all $x \in \Omega$, $h_t(x)= x$ implies that $t=e$, where $e$ is the group unit.
\end{definition}

It is well known that the category of partial actions on a Hausdorff locally compact space $\X$ is equivalent to the category of partial actions on the C*-algebra of the continuous functions vanishing at infinity $\CC_0(\X)$ (see proposition 1.5 of \cite{Abadie} for example). In our work we will use that given a partial action $(\{\X_t\}_{t \in G}, \{h_t\}_{t \in G})$ of $\G$ on $\X$, then $(\{\CC_0(\X_t)\}_{t \in G}, \{\alpha_t\}_{t \in G})$, where $\alpha_t :\CC_0(\X_{t^{-1}})  \rightarrow  \CC_0(\X_{t})$ is given by $\alpha_t(f)=f \circ h_t^{-1}$, is a partial action of $\G$ on $\CC_0(\X)$. Next we recall the definition of the partial crossed product.

\begin{definition}
Let $\theta = (\{ D_{t} \}_{t \in \G} , \{ \alpha_{t} \}_{t \in \G} )$ be a partial action of $\G$ on the C*-algebra $\A$. Then the partial crossed product of $\A$ by $\G$, denoted by $A \rtimes_\alpha G$, is the enveloping $C^*$-algebra of $\mathcal{L}$, where $\mathcal{L}$ is the normed *-algebra of all finite sums of the form $$\mathcal{L}=
\left\lbrace \sum_{t \in G}a_t\delta_t \,\, : \,\, a_t \in
\mathcal{D}_t \,\,\, \right\rbrace  \subseteq l_1(G, A)$$ where, for $a=(a_t)_{t\in
G} \in \mathcal{L}$ and $b=(b_t)_{t\in G} \in \mathcal{L}$, multiplication, involution and norm are defined by  
$(a \ast
b)_{\gamma}= \sum_{t \in
G}\alpha_t(\alpha_{t^{-1}}(a_t)b_{t^{-1}\gamma})$, 
$(a^*)_{\gamma}=\alpha_{\gamma}(a_{\gamma^{-1}}^*)$ and $\Vert a \Vert = \sum_{t \in G} \Vert a_t\Vert.$
\end{definition}

Another key definition we need to recall is that of an étale equivalence relation. In the language of \cite{Renault} this is an r-discrete groupoid with counting measure as a Haar system. In our context an étale equivalence relation $R\subseteq \X \times \X$, where $\X$ is  locally compact Hausdorff, is an equivalence relation that can be equipped with two maps, called range and source, defined by $r\left(x,y \right)=x$ and $s\left(x,y \right)=y $ and such that $R$ is $\sigma$-compact, $\Delta = \{(x,x) \in R : x\in \X \}$ is an open subset of $R$ and for all $\left(x,y \right) \in R$, there exists a neighborhood $\U$ of $\left(x,y\right)$ in $R$, such that $r$ restricted to $\U$ and $s$ restricted to $\U$ are homeomorphisms from $\U$ onto open subsets of  $\X $, see also \cite{Putnam3}.

Finally, given an étale equivalence relation $R$ over a locally compact set $\X$, the groupoid algebra associated to it is obtained as the completion, over a certain norm (see \cite{Renault} or \cite{Gonc}), of the *-algebra of the continuous functions with compact support in $R$, $\CC_c(R)$, where the *-algebra operations are defined, for $f, g \in \CC_c(R)$, as
$f \ast g (x,z) =\displaystyle \sum_{y \in [x]} f(x,y)g(y,z)$ and $ 
f^*(x,y)= \overline{f(x,y)}$ (where $[x]$ denotes the equivalence class of $x$).

\section{C*-algebra level}

Let $(\{ X_t\}_{t \in \G}, \{h_{t \in \G}\}) $ be a free partial action of a countable group $\G$ on a $\sigma$-compact Hausdorff space $\X$ such that $X_t$ is $\sigma$-compact for every $t$. In this section, we prove that the partial crossed product $\CC(X)\rtimes _{\alpha} G$ associated to the corresponding partial action $(\{\CC(\X_t)\}_{t \in G}, \{\alpha_t\}_{t \in G})$ (as defined in the introduction)  is isomorphic to the full groupoid C*-algebra $\CC^*(R)$, where $R$ is defined below.

\begin{definition}\label{relacao} We say that $x$ is equivalent to $y$, $x \sim y $, if there exists $t \in G$ such that $x \in X_{t^{-1}}$ and $h _t(x)=y$. 
\end{definition} 

\begin{obs}
Notice that the elements of $R$ are of the form $(x, h_t(x))$ and, since the action is free, for each $(x,
y)\in R$ there exists one and only one $t \in G$ such that $(x, y)= (x, h_t(x))$.
\end{obs}

Of course we are still missing the key ingredient before we can proceed. That is the topology of $R$. We stress out that this is not the topology inherited from the product topology of $\X \times \X$. Instead, to obtain an étale equivalence relation, we give $R$ the inherited topology from $\X \times \G$ (with the product topology) via the map $ (x, h_t(x)) \in R \longmapsto (x, t) \in \X \times \G $. Notice that since the action is free, this map is injective. So, a sequence $\{(x_n, h_{t_n}(x) \}$ converges to $(x, h_t(x)) \in R$ iff $\{ x_n \}$ converges to $x$ in $\X$ and $t_n$ are eventually all equal to $t$. With this topology, we can now prove that $R$ is étale.

\begin{proposition} $R$ is an étale equivalence relation.
\end{proposition}
\demo

Before we show that $R$ is étale we should prove that $R$ is an equivalence relation. It is straightforward to check that $R$ is symmetric and reflexive. We show that $R$ is transitive below.

If $x \sim y $ and $y \sim z$ then there exist $t, s \in G$ such that $x \in X_{t^{-1}}, y \in X_t\cap X_{s^{-1}}, z \in X_s , h_t(x)=y$ and $ h_s(y)=z $. This implies that $x \in h_{t}^{-1}(X_t \cap X_{s^{-1}})\subseteq X_{(st)^{-1}}$, $z \in h_s( X_{s^{-1}}\cap X_t)=X_s \cap X_{st}$ and $z=h_s(y)= h_s(h_t(x)) =h_{st}(x) $. Taking $r=st$ we have that $x \in X_{r^{-1}}$, $z \in X_r$ and $h_r(x)=z$. We conclude that $x \sim z$.

Now we prove that $R$ is étale. By hypothesis, we have that $X_t$ is $\sigma$-compact for all $t \in G$. Then, for each $t$, there exists a countable family of compact subsets of $X_t$, $\{K_n^t\}_{n \in \N}$, s.t. $X_t= \bigcup_{n \in \N}K_n^t$. We conclude that the sets $$U_{t, n}:= \{(x, h_t(x)) \,\, : \,\, x \in K_n^t\}$$
are compact and  $R= \bigcup_{t \in G, n \in \N}U_{t, n}$ is
$\sigma$-compact. That the diagonal $\bigtriangleup = \{(x, h_e(x))
\,\, : \,\, x \in X \} = U_e$ is open in $R$ is clear. Finally, if $(x, h_t(x)) \in R$ then $U_t:= \{(z, h_t(z)) \,\, : \,\, z \in X_{t^{-1}}\}$ is an open neighbourhood of $(x, h_t(x))$ and the range $r: U_t \rightarrow X_{t^{-1}}$ and source $s: U_t \rightarrow X_t$ maps are homeomorphisms.

\cqd

\begin{obs} Notice that, since the partial action $(\{ X_t\}_{t \in \G}, \{h_{t \in \G}\}) $ is free, the map $(x,y) \mapsto (y,t,x)$, where $t$ is the unique element of $\G$ such that $h_t(x)=y$,  is an isomorphism from $R$ onto $\mathcal{G}$, where $\mathcal{G}$ is the groupoid constructed in \cite{Abadie}.
\end{obs}

We can now consider $\CC^*(R)$ and show that it is isomorphic to $\CC(X)\rtimes _{\alpha} G$. Our proof relies on the theory of core sub-algebras, so that, by proposition 3.4 in \cite{ExGiGo}, it is enough to show that $\mathcal{L}$ and $\CC_c(R)$ are dense core sub-algebras of  $\CC(X)\rtimes _{\alpha} G$ and  $\CC^*(R)$, respectively, and that they are isomorphic *-algebras. For reader´s convenience we recall the definition of a core sub-algebra below.

\begin{definition}Let $A$ be a C*-algebra and let $B\subseteq A$ be a (not necessarily
closed) *-subalgebra.  We say that $B$ is a core
subalgebra of $A$ when every representation \footnote{By a representation of a
*-algebra $B$ we mean a multiplicative, *-preserving, linear map
$\pi:B\to{\mathcal B}(H)$, where $H$ is a Hilbert space.}
  of $B$ is continuous relative to the norm induced from $A$.
\end{definition}

As a consequence of our next result we obtain that $\mathcal{L}$ is a core sub-algebra of $\CC(X)\rtimes _{\alpha} G$. For the analogous result concerning $\CC_c(R)$ we refer the reader to \cite{ExGiGo} and \cite{Renault}.

\begin{proposition} Let $(\{ D_t\}_{t \in G}, \{\alpha_{t \in G}\}) $ be a partial action of a discrete group $\G$ over the C*-algebra $A$. Then 
$\mathcal{L}$ is a core sub-algebra of $A \rtimes _{\alpha} G$.
\end{proposition}
\demo
Let $\lambda$ be a representation of $(\mathcal{L}, |\Vert \cdot \Vert|)$. We have to show that $\Vert \lambda(a) \Vert \leq |\Vert a \Vert|$ for all $a \in \mathcal{L}$, where $$|\Vert a \Vert| = \displaystyle\sup_{\pi}\{\Vert \pi(a) \Vert \, : \, \pi \, \mbox{is a representation s.t.} \,\,\Vert \pi(a) \Vert \leq \Vert a \Vert_{l_1} \}.$$ 
Notice that $A_0= \{ a\delta_0 \in \mathcal{L} \,\, : \,\, a \in A\}$, equipped with the  $\mathcal{L}$ operations, is a $C^*$-algebra (isomorphic to A) and so $\Vert \lambda(a\delta_0 ) \Vert \leq \Vert a \delta_0 \Vert_{l_1} \,\, \mbox{for all}  \,\,\, a\delta_0 \in \mathcal{L} .$
It follows that,
\begin{align*}
\Vert \lambda (a_t \delta_t) \Vert^2 & = \Vert  \lambda(a_t \delta_t )\lambda^*(a_t \delta_t)\Vert = \Vert \lambda((a_t\delta_t)(a_t\delta_t)^*)\Vert \\ & = \Vert \lambda((a_t\delta_t)\alpha_{t^{-1}}(a_t^*)\delta_{t^{-1}}) \Vert =   \Vert \lambda(a_ta_t^*\delta_0)\Vert \\ & \leq \Vert a_ta_t^*\delta_0 \Vert = \Vert a_ta_t^* \Vert  = \Vert a_t \Vert^2= \Vert a_t\delta_t \Vert^2,
\end{align*}
and hence $\Vert \lambda (a_t \delta_t) \Vert \leq \Vert a_t\delta_t \Vert_{l_1} \,\, \mbox{for all} \,\, a_t\delta_t \in \mathcal{L}.$ \\
Now let $a= \displaystyle{\sum_{t \in G}}a_t \delta_t \in \mathcal{L}$. Then
\begin{align*}
\Vert \lambda( a) \Vert & = \Vert \lambda (\displaystyle{\sum_{t \in G}}a_t \delta_t) \Vert \leq \displaystyle{\sum_{t \in G}}\Vert \lambda(a_t\delta_t) \Vert  \leq \displaystyle{\sum_{t \in G}}\Vert a_t\delta_t \Vert \\ & = \displaystyle{\sum_{t \in G}}\Vert a_t \Vert = \Vert  a \Vert _{l_1}.
\end{align*}
and so  $\Vert \lambda( a) \Vert  \leq \Vert  a \Vert _{l_1},  \,\, \mbox{for all} \,\, a \in \mathcal{L}.$ 

We conclude that $\lambda$ is one of the representations over which the sup in the definition of $|\Vert a \Vert|$ is taken and hence $\Vert \lambda(a)\Vert \leq  |\Vert a \Vert|$ for all $a \in \mathcal{L}$.

\cqd

\begin{obs} The idea presented above can also be used to show that
$l_1(\G,\A)$ is a core sub-algebra of the full crossed product $A \rtimes _{\alpha} G$ associated to an action of $\G$ on $\A$.
\end{obs}

Now that we have estabilished that $\mathcal{L}$ is a core sub-algebra of  $\CC(X)\rtimes _{\alpha} G$ all we are left to do, by proposition 3.4 in \cite{ExGiGo}, is to show that  $\mathcal{L}$ is isomorphic to $\CC_c(R)$.

\begin{theorem} Let  $(\{ X_t\}_{t \in \G}, \{h_{t \in \G}\}) $ be a free partial action of a countable group $\G$ over a locally compact Hausdorff space $\X$ such that $X_t$ is $\sigma$-compact for every $t$, $(\{\CC(\X_t)\}_{t \in G}, \{\alpha_t\}_{t \in G})$ be the corresponding partial action (as defined in the introduction) and $R$ be the equivalence relation defined above. Then $\mathcal{L}$ and $\CC_c(R)$ are isomorphic as *-algebras.
\end{theorem}
\demo

To define a *-homomorphism $\rho : \mathcal{L} \longrightarrow   C_c(R)$ we begin defining it at elements of the form $f_t \delta_t$ ($f_t \in \CC_0(X_t)$) and then we extend it linearly to $\mathcal{L}$. More precisely, for  $f_t \delta_t \in \mathcal{L}$ and $(x,h_s(x)) \in R$, let $$\tilde{\rho} (f_t \delta_t ) (x,h_s(x)) =  \left\lbrace \begin{array}{ll}
                                    f_t(x)  &  \mbox{if} \,\,\,  s=t^{-1} \\
                                    0  &  \mbox{otherwise}
                                 \end{array}\right.,$$ and denote the linear extension of $\tilde{\rho}$ to $\mathcal{L}$ by $\rho$.

Notice that since $\mathcal{L}$ consists of finite sums the set $\{s\in \G: \rho (\sum f_t \delta_t ) (x,h_s(x)) \neq 0 \}$ is finite and hence $\rho$ is well defined (that is, $\rho(\sum f_t \delta_t)$ is a continuous function with compact support). 

Next we check that $\rho$ is *-multiplicative. By linearity, it is enough to check this for elements of the form $f_t \delta_t \in  \mathcal{L}$.
So, let $f_t \delta_t, g_s \delta_s \in  \mathcal{L}$ and $(x,h_r(x))\in R$. Then,

\begin{align*}
  \rho(f_t\delta_t\ast g_s\delta_s)(x,h_r(x)) & = \rho(\alpha_t(\alpha_{t^{-1}}(f_t)g_s) \delta_{ts}) (x,h_r(x))\\
                   & = \left\lbrace \begin{array}{ll}
                                    f_t(x)g_s(h_{t^{-1}}(x))  &  \mbox{if} \,\,\,  r=ts \\
                                    0  &  \mbox{otherwise}
                                 \end{array}\right.\\
							& = \left\lbrace \begin{array}{ll}
                                    f_t(x)\rho(g_s\delta_s)(h_{t^{-1}}(x), h_r(x))  &  \mbox{if} \,\,\,  r=ts \\
                                    0  &  \mbox{otherwise}
                                 \end{array}\right.\\
 						& = \displaystyle \sum_{u: x\in X_{-u}} \rho(f_t \delta_t)(x,h_u(x)) \rho(g_s \delta_s) (h_u(x), h_r(x)) \\
 						& =\rho(f_t \delta_t) * \rho (g_s \delta_s) (x, h_r(x)),
\end{align*} and
\begin{align*}
\rho((f_g\delta_g)^*)(x,h_r(x)) & = \rho(\alpha_{g^{-1}}(f_g^*)\delta_{g^{-1}}) (x, h_r(x)) \\
             &          = \left\lbrace \begin{array}{ll}
                                  \alpha_{g^{-1}}(f_g^*)(x)  &  \mbox{if} \,\,\,  r=g \\
                                    0  &  \mbox{otherwise}
                                 \end{array}\right.\\
      &                  = \left\lbrace \begin{array}{ll}
                                 \overline{f_g(h_r(x))}  &  \mbox{if} \,\,\,  r=g \\
                                    0  &  \mbox{otherwise}
                                 \end{array}\right.\\
   & = \overline{\rho(f_g \delta_g) (h_r(x),x)} = (\rho(f_g\delta_g))^*(x,h_r(x)).
\end{align*}

So $\rho$ is a *-homomorphism. Finally, notice that $\rho$ is also a bijection, since it has an inverse $\rho^{-1}: C_c(R) \longrightarrow \mathcal{L}  $  given by $\rho^{-1}(f) = \sum f_t \delta_t$, where $f_t(x)=  \left\lbrace \begin{array}{ll}
                                 f(x,h_{t^{-1}}(x))  &  \mbox{if} \,\,\,  (x, h_{t^{-1}}(x)) \in R \\
                                    0  &  \mbox{otherwise}
                                 \end{array}\right..$

\cqd

\begin{obs} Notice that since an action of a group on an algebra can be seen as a partial action the above result is also valid in the context of actions of countable groups over compact spaces.
\end{obs}

\begin{corolario} $\CC(X)\rtimes _{\alpha} G$ is isomorphic to $\CC^*(R)$.
\end{corolario}
\demo
Follows from the above theorem and proposition 3.4 in \cite{ExGiGo}.

\cqd

\section{The purely algebraic setting}

In this section we generalize to the purely algebraic setting the results of the previous section and the correspondence between partial dynamical systems and partial C*-dynamical systems.

Let $\K$ be a field and $\X$ a set. By an algebra we mean an associative $\K$-algebra, not necessarily unital. Let $\Fo$ denote the algebra of all functions $f: \X \rightarrow \K$ that vanish eventually, that is $f\in \Fo$ iff $f(x)=0$ for all but a finite number of $x \in \X$ ($f$ has finite support), equipped with pointwise operations. Notice that we can see $\Fo$ as the direct sum of $\K$ over $\X$, but we will keep the function notation due to it resemblance to the C*-setting. Our first goal is to show that there exists a bijective correspondence between $X$ and the set of all non zero homomorphisms from $\Fo$ to $\K$, where each $x \in X$ is taken to the homomophism $\epsilon_x : \Fo  \rightarrow  \K $, given by $\epsilon_x(f)=f(x)$, which we call evaluation at $x$.

\begin{proposition}\label{biunivoca}
There exists a bijective correspondence between a set $\X$ and $\widehat{\Fo}$ given by:
$ \begin{array}{rcl}
\X & \rightarrow & \widehat{\Fo} \\
x & \mapsto & \epsilon_x
\end{array},$ where $\widehat{\Fo}$ denotes the set of all non zero homomorphisms from $\Fo$ to $\K$.
\end{proposition}

\demo

Given $x\in \X$, let $\delta_x$ denote the characteristic function of the set $\{x\}$. Then $\{\delta_x\}_{x\in \X}$ is a $\K$-basis of $\Fo$ and $f= \sum_{x\in \X} f(x) \delta_x$, for all $f \in \Fo$.

First we will prove that there is a bijective correspondece between $\K$-linear maps from $\Fo$ to $\K$ and functions from $X$ to $\K$. For this notice that, if $\phi: \Fo \rightarrow \K$ is a $\K$ linear map and $f\in \Fo$, then $$\phi(f)= \sum_{x\in \X}f(x) \phi (\delta_x) = \sum_{x\in \X}f(x) \varphi(x)=\sum_{x\in \X}\varphi(x)\epsilon_x(f),$$ where $\varphi: X \rightarrow \K$ is defined by $\varphi(x)= \phi (\delta_x), \  \ \forall x \in \X$. Conversely, every function $\varphi: \X \rightarrow \K$ defines a $\K $-linear map $\phi: \Fo \rightarrow \K$ via $\phi(f) = \sum_{x \in \X} \varphi(x) \epsilon_x(f)$, and it is clear that the correspondence $\phi \leftrightarrow \varphi$ is a bijection.

Now, notice that every evaluation $\epsilon_x$ is a nonzero homomorphism. Suppose conversely that $\phi: \Fo \rightarrow \K$ is a nonzero homomorphism. Then we have that $\varphi(x)\varphi(y)= \phi(\delta_x) \phi (\delta_y) = \phi (\delta_x \delta_y) = \delta_{x,y} \varphi(x)$. Then $\varphi(x) \varphi(y) = 0 $ if $x\neq y$ and $\varphi(x)^2 = \varphi(x)$ for all $x\in \X$. Hence there is exactly one $x\in X$ such that $\varphi(x) \neq 0$, so $\varphi = \epsilon_x$.

\cqd

\par Our next goal is to show that there is a biunivocal correspondence between partial actions $(\{X_t\}_{t \in G}\{h_t\}_{t \in G})$ of a group $G$ in a set $X$ and partial actions $(\{D_t  \}_{t \in G}, \{\alpha_t\}_{t \in G})$ of $G$ in $\Fo$. Before we do this we need a few results.


\begin{proposition}\label{isomorfismo}
Let $X$ and $Y$ be sets and $h: X \rightarrow Y$ a bijection. Then the map $\psi_h: \mathcal{F}_0(Y)   \rightarrow  \Fo$, defined by $\psi_h(f)= f\circ h$, is an algebra isomorphism.
\end{proposition}
\demo
The proof of this propostition is  straighforward.
\cqd

\begin{proposition}\label{bijecao}
If $\gamma: \mathcal{F}_0(Y) \rightarrow \Fo$ is an isomorphism then there exists a unique bijection, $h: X \rightarrow Y$, such that $\gamma = \psi_h$, where $\psi_h$ is as in the previous proposition.
\end{proposition}

\demo

Let $\gamma: \mathcal{F}_0(Y) \rightarrow \Fo$ be an isomorphism. We need to define $h:X \rightarrow Y$. For this, notice that for all $x \in X$, ${\epsilon_x}\circ\gamma$ is a homomorphism from $\mathcal{F}_0(Y)$ to $\K$. By proposition \ref{biunivoca} all such homomorphisms are evaluations and hence there exists $y \in Y$ such that ${\epsilon_x}\circ\gamma=\epsilon_y$. Define $h(x):=y$. 

Notice that indeed $\gamma= \psi_h$, since for all $f \in \mathcal{F}_0(Y) $ and $x \in X$, we have that $\psi_h(f)(x)=f(h(x))=f(y)=\epsilon_y(f)={\epsilon_x}\circ\gamma(f)=\epsilon_x(\gamma(f)) =\gamma(f)(x)$. Also, $h$ is bijective, since we can define its inverse, $l: Y \rightarrow X $, in the following way: Consider $\gamma^{-1}: \Fo \rightarrow \mathcal{F}_0(Y)$. Given $y \in Y$, ${\epsilon_y}\circ\gamma^{-1}$ is a homomorphism from $\Fo$ to $\K$ and so there exists $x \in X$ such that ${\epsilon_y}\circ\gamma^{-1}=\epsilon_x$. Define $l(y):=x$. Then, for $x \in X$, $l\circ h(x)=l(y)=x_0$, where $y$ is such that $\epsilon_x \circ \gamma= \epsilon_y$ and $x_0$ is such that $\epsilon_y\circ\gamma^{-1}= \epsilon_{x_0}$. So,
$x_0$ is such that $\epsilon_{x_0}=\epsilon_y \circ \gamma^{-1}=\epsilon_x \circ \gamma \circ \gamma^{-1}= \epsilon_x$ and hence $x_0=x$ and $l \circ h = Id$. Analogously we can check tat $h\circ l= Id$ and hence $l = h^{-1}$.

Finally we show that there is a unique $h$ such that $\gamma= \psi_h$. For this, suppose that there exists bijections $h_1$ and $h_2$ such that $\gamma(f)=\psi_{h_1}(f)=\psi_{h_2}(f)$ for all $f \in \mathcal{F}_0(Y)$. We then have that $f \circ h_1 (x) = f \circ h_2(x)$ for all $f \in \mathcal{F}_0(Y),\, x \in X$, what implies that $h_1(x)= h_2(x)$ for all $x \in X$, since if there exists $x \in X$ such that $h_1(x)\neq h_2(x)$ then, for $f=\delta_{h_1(x)}$, we have that $f(h_1(x))\neq f(h_2(x))$. \cqd 

\begin{proposition}\label{ideals} There is a biunivocal correspondence between non zero ideals of $\Fo$ and non empty subsets of $X$. 
\end{proposition}
\demo
Let $I$ be a non zero ideal in $\Fo$. Define $A$ as the set of all elements of $X$ such that there exists a function $f\in I$ such that $\delta_x \cdot f \neq 0$. Then $I=\Fa$, where $\Fa$ is included in $\Fo$, that is, $\Fa=\{f\in \Fo: \ f(x)=0, \ \forall \ x\notin A \}$.
\cqd\\

From the above propositions, we conclude that there is a bijective homomorphism $\psi$ between the set of all bijections from $X$ to $Y$, which we denote by  $Bij(X, Y)$, and the set of all isomorphisms from $\Fo$ to $\mathcal{F}_0(Y)$, which we denote by $Iso( \Fo, \mathcal{F}_0(Y))$, given by
$$\begin{array}{rcl}
\psi: Bij(X, Y) & \rightarrow & Iso( \Fo, \mathcal{F}_0(Y))\\
              h & \mapsto     & \psi_{h^{-1}}
\end{array}.$$ 

We are now ready to prove the correspondence between partial actions on a set $X$ and partial actions on $\Fo$.

\begin{proposition}\label{actionarises}
Let $\theta=(\{X_t\}_{t \in G}, \{h_t\}_{t \in G})$ be a partial action of a group $G$ in a set $X$ and let $D_t=\{f\in \Fo: \ f(x)=0 \ \forall \ x\notin X_t \}$, that is, $D_t = \mathcal{F}_0(X_t)$. Define $\alpha_t : D_{t^{-1}} \rightarrow D_t$ by $\alpha_t(f)=f\circ h_{t^{-1}}$. Then $\alpha=(\{D_t\}_{t \in G}, \{\alpha_t\}_{t \in G})$ is a partial action of $G$ in $\Fo$ and we say that $\alpha$ arises from $\theta$.  
\end{proposition} 

\demo It is clear that, for each $t\in G$, $D_t$ is an ideal of $\Fo$ and, by proposition \ref{isomorfismo}, $\alpha_t$ is bijective. We show below that $(\{D_t\}_{t \in G}, \{\alpha_t\}_{t \in G})$ satisfies the other axioms of the definition of a partial action. 

\begin{itemize}
\item $D_e = \Fo$ and $\alpha_e = Id$.\\
Since $\theta=(\{X_t\}_{t \in G}, \{h_t\}_{t \in G})$ is a partial action we have that $X_e=X$ and $h_e=Id$ what readily implies that $D_e=\Fo$ and $\alpha_e = Id$.

\item $\alpha_t(D_{t^{-1}}\cap D_s) = D_t \cap D_{ts}$.\\
First we prove that $\alpha_t(\mathcal{F}_0(X_{t^{-1}}\cap X_s))= \mathcal{F}_0(h_t(X_{t^{-1}}\cap X_s))$. For this, let $y=f\circ h_{t^{-1}}$, for some $f\in \mathcal{F}_0(X_{t^{-1}}\cap X_s)$. Now, if $x\notin h_t(X_{t^{-1}}\cap X_s)$ then $h_{t^-1}(x) \notin X_{t^{-1}}\cap X_s$ and hence $f\circ h_{t^{-1}} (x)=0$ and $y \in \mathcal{F}_0(h_t(X_{t^{-1}}\cap X_s))$. On the other hand, if $f\in \mathcal{F}_0(h_t(X_{t^{-1}}\cap X_s))$ then $f\circ h_t \in \mathcal{F}_0(X_{t^{-1}}\cap X_s)$ and $\alpha_t(f\circ h_t) = f$ as desired.

We can now prove the desired partial action aximom below:

$\alpha_t(D_{t^{-1}}\cap D_s)= \alpha_t(\mathcal{F}_0(X_{t^{-1}})\cap \mathcal{F}_0(X_s))= \alpha_t(\mathcal{F}_0(X_{t^{-1}}\cap X_s))= \mathcal{F}_0(h_t(X_{t^{-1}}\cap X_s))=\mathcal{F}_0(X_t\cap X_{ts}) = \mathcal{F}_0(X_t)\cap \mathcal{F}_0(X_{ts})=D_t \cap D_{ts}$.

\item $\alpha_t(\alpha_s(f)) = \alpha_{ts}(f)$ for all $f \in D_{s^-1}\cap D_{s^{-1}t^{-1}}$.\\
Let $f\in D_{s^{-1}}\cap D_{s^{-1} t^{-1}} = \mathcal{F}_0(X_{s^{-1}} \cap \X_{s^{-1} t^{-1}})$ and $ x \in X_{s^{-1}} \cap \X_{s^{-1} t^{-1}} $. It follows that:

$\alpha_t(\alpha_s(f))(x)=\alpha_t(f\circ h_{s^{-1}})(x)=f\circ h_{s^{-1}}\circ h_{t^{-1}}(x)=f(h_{s^{-1}}( h_{t^{-1}}(x)))=f(h_{s^{-1}t^{-1}})(x)= \alpha_{ts}(f)(x)$.
\end{itemize}
\cqd \\

\begin{proposition} 
If $\alpha=(\{D_t\}_{t \in G}, \{\alpha_t\}_{t \in G})$ is a partial action of  $G$ in $\Fo$ then there exists a partial action $\theta=(\{X_t\}_{t \in G}, \{h_t\}_{t \in G})$, of $G$ in a set $X$, such that $\alpha$ arises from $\theta$.
\end{proposition}

\demo
Let $\alpha=(\{D_t\}_{t \in G}, \{\alpha_t\}_{t \in G})$ be a partial action of $G$ in $\Fo$. By proposition \ref{ideals} we have that each ideal $D_t$ is of the form $\mathcal{F}_0(X_{t})$, for some subset $X_t$ of $X$. Now, for each isomorphism $\alpha_t :  \mathcal{F}_0(X_{t^{-1}}) \rightarrow \mathcal{F}_0(X_t)$, we let $h_{t^{-1}}$ be the unique bijection from $X_t$ to $X_{t^{-1}}$ such that $\alpha_t = \psi_{h_{t^{-1}}}$  (that is, $\alpha_t(f)=f\circ h_{t^{-1}}$ as described in proposition \ref{bijecao}). 
This way we define a partial action $\theta= (\{X_t\}_{t \in G}, \{h_t\}_{t \in G})$ such that $\alpha$ arises from $\theta$. To finish the prove we need to we show that $\theta $ is indeed a partial action. 

It is straightforward to check that $\theta $ satisfies the first axiom in the definition of a partial action. To verify that $h_t(X_{t^{-1}}\cap X_s) = X_t \cap X_{ts}$ notice that
$\alpha_t(\mathcal{F}_0(X_{t^{-1}})\cap\mathcal{F}_0( X_s))= \mathcal{F}_0(X_t)\cap \mathcal{F}_0(X_{ts})=\mathcal{F}_0(X_t \cap X_{ts})$ and since the left side on this last equality is equal to $\alpha_t(\mathcal{F}_0(X_{t^{-1}}\cap X_s)) = \mathcal{F}_0(h_t(X_{t^{-1}}\cap X_s)$, we obtain the desired equality.
Finally we prove that $h_t(h_s(x)) = h_{ts}(x)$ for all $x \in X_{s^-1}\cap X_{s^{-1}t^{-1}}$. Since $\alpha$ is a partial action, we have that $\alpha_t(\alpha_s(f)) = \alpha_{ts}(f)$ for all $f\in \mathcal{F}_0 (X_{s^-1}\cap X_{s^{-1}t^{-1}})$. But this implies that $\psi_{h_{ts}}(f)=  \psi_{h_t \circ h_s}(f)$ for all $f\in \mathcal{F}_0 (X_{s^-1}\cap X_{s^{-1}t^{-1}})$ and hence $h_{ts}=h_t \circ h_s$ for all $x \in X_{s^-1}\cap X_{s^{-1}t^{-1}}$ as desired.
\cqd \\

 We now focus on realizing the partial skew group ring $\Fo \rtimes_{\alpha} G$, associated to a partial action $\alpha$ of a group $G$ in $\Fo$, as an algebra from an equivalence relation. 

For reader´s convenience we recall the definition of a partial skew group ring, as defined in \cite{Ex}, below.

\begin{definition}
Let $\alpha$ be a partial action of the group $G$ in the algebra $A$. The partial skew group ring, $A\rtimes _{\alpha} G$, associated to $\alpha$ is
defined as the set of all finite formal sums $\sum_{t\in G} a_t \delta_t$,  where, for all $t\in G$, $a_t \in D_t$ and $\delta_t$ are symbols. Addition is
defined in the usual way and multiplication is determined by
 $(a_t\delta_t)(b_s\delta_s)=\alpha_t(\alpha_{t^{-1}}(a_t)b_s)\delta_{ts}$.
\end{definition}

Recall that for every partial action in $\Fo$ there is a partial action in $X$ associated. The set level is what we need to define the corresponding  equivalence relation. So, let $(\{X_t\}_{t \in G}, \{h_t\}_{t \in G})$ be a partial action of $G$ in $X$. We define the equivalence relation $R$ in $X$ as
  $$R=\{(x, h_t(x))\in X \times X\, : \, x \in X_{t^-1}, t \in G\},$$
and equip
$\Fr = \{f: R \rightarrow \K \, : \, f \, \mbox{eventually vanishes}\}$
with operations defined by:

$(kf+g)(x,h_t(x))=kf(x,h_t(x))+g(x, h_t(x))$, 
$$f \ast g (x,h_t(x))= \sum_{s \in G}f(x,h_s(x))g(h_s(x),h_t(x)),$$
for all $f, g \in \Fr$,  $k \in \K$ and $(x, h_t(x)) \in R$

\begin{theorem}\label{oteorema}
If $\theta=(\{X_t\}_{t \in G}, \{h_t\}_{t \in G})$ is a free partial action of a group $G$ in a set $X$, and $\alpha = (\{D_t\}_{t \in G}, \{\alpha_t\}_{t \in G})$ is the associated partial action of $G$ in $\Fo$, then the algebras $\Fr$ and $\Fo \rtimes_{\alpha}G$ are isomorphic.
\end{theorem}

\demo
To prove the theorem we will define an isomorphism $\Gamma$ from $\Fo \rtimes_{\alpha} G$ to $\Fr$. For elements of the form $f\delta_t \in \Fo \rtimes_{\alpha} G$ we let $\Gamma(f\delta_t)=\widetilde{f}$, where
$$\widetilde{f}(x,y)=\left\lbrace  \begin{array}{ll}
                        f(x) & \mbox{if} \,\,\, y=h_{t^{-1}}(x)\\
                        0      & \mbox{otherwise}
                       \end{array},\right. $$
and we define $\Gamma$ in $\Fo \rtimes_{\alpha} G$ by extending it linearly. 

Notice that $\Gamma$ is well defined since if $\{x_1^t, \cdots, x_{n_t}^t\}$ is the finite support of $f_t \in \Fo$ then $\Gamma(f_t\delta_t)$ has finite support given by the set $\{(x_1^t, h_{t^{-1}}(x_1^t)), \cdots , (x_{n_t}^t, h_{t^{-1}}(x_{n_t}^t))\}$. 

Now we show that $\Gamma$ is an isomorphism of $\K$-algebras. For this it is enough to show that $\Gamma$ is bijective and preserves the product, that is,  $\Gamma(f_t \delta_t \ast f_g\delta_g)= \Gamma(f_t \delta_t) \ast \Gamma(f_g \delta_g)$, for all $f_t\delta_t, f_g\delta_g \in \Fo \rtimes_{\alpha} G$. Notice that,
\begin{align*}
\Gamma(f_t\delta_t \ast f_s \delta_s)(x,y) & = \Gamma(\alpha_t(\alpha_{t^{-1}}(f_t) f_s) \delta_{ts})(x,y)\\
 & = \left\lbrace \begin{array}{ll}
       \alpha_t(\alpha_{t^{-1}}(f_t) f_s)(x), & \mbox{if } \, y=h_{(ts)^{-1}}(x)\\
        0                     & \mbox{otherwise}
      \end{array}\right.\\
 & = \left\lbrace \begin{array}{ll}
       f_t(x)(f_s)(h_{t^{-1}}(x)), & \mbox{if } \, y=h_{(ts)^{-1}}(x)\\
        0                     & \mbox{otherwise}
      \end{array}\right.\
\end{align*}
And, on the other hand,
\begin{align*}
\Gamma(f_t\delta_t) \phi(f_s \delta_s)(x,y) & = \sum_{r \in G} \Gamma(f_t \delta_t)(x, h_r(x))\Gamma(f_s \delta_s)(h_r(x),y)\\
 & = f_t(x)\Gamma(f_s\delta_s)(h_{t^{-1}}(x), y)\\
 & = \left\lbrace \begin{array}{ll}
       f_t(x)(f_s)(h_{t^{-1}}(x)), & \mbox{if } \, y=h_{s^{-1}}(h_{t^{-1}}(x))\\
        0                     & \mbox{otherwise}
      \end{array}\right.\\
 & = \left\lbrace \begin{array}{ll}
       f_t(x)(f_s)(h_{t^{-1}}(x)), & \mbox{if} \, y=h_{(ts)^{-1}}(x)\\
        0                     & \mbox{otherwise}
      \end{array}\right.
\end{align*}
We conclude that, $\Gamma(f_t\delta_t \ast f_s \delta_s) =\Gamma(f_t\delta_t) \Gamma(f_s \delta_s)$. 

Next we prove injectivity. For this, supose that $ \Gamma\left( \sum_{t \in G}f_t \delta_t\right)  = 0.$ Notice that, since $\theta$ is a free, if $\Gamma (f_t \delta_t)(x,y)\neq 0$ then $\Gamma (f_s \delta_s)(x,y)= 0$ for all $s\neq t$. So, $  \sum_{t \in G}\Gamma(f_t \delta_t)(x,h_s(x))=0$, for all $x \in X_{s^{-1}},  s \in G$ implies that $ f_{s^{-1}}(x)=0$ for all $x \in X_{{s^-1}},  s \in G$. We infer that $ f_t \equiv 0 $ for all $ t \in G $ and hence $ \sum_{t \in G} f_t \delta_t = 0 $ as desired.

Finally we prove that $\Gamma$ is surjective: for $f \in \Fr$ we define $f_t(x): = f(x, h_{t^{-1}}(x))$, for all $x \in X_t$. 
Notice that $\sum_{i=1}^n f_{t_i}\delta_{t_i} \in \Fo \rtimes_{\alpha} G$ and
$\Gamma\left( \sum_{i=1}^n f_{t_i}\delta_{t_i}\right) = f$ so that $\Gamma$ is surjective.

By the exposed above we conclude that $\Fo \rtimes_{\alpha}G$ and $\Fr$ are isomorphic as $\K$-algebras. \cqd

\begin{obs} Partial skew group rings similar to the ones considered in the above theorem have appeared before, see \cite{GR1}. 
\end{obs}

One important issue that has to be handled when dealing with partial skew group rings concerns the associativity question, as partial skew group rings may or may not be associative, see \cite{Ex}. In the case of $\Fo \rtimes_{\alpha}G$ the answer is affirmative, that is, $\Fo \rtimes_{\alpha}G$ is always associative, since all the ideals of $\Fo$ are idempotent (see \cite{Ex,BR} for associativity criteria). Below we give a quick proof that $\Fo \rtimes_{\alpha}G$ is always associative, one that does not depend in the theory developed in \cite{Ex, BR}. 

\begin{proposition} $\Fr$ is an associative $\K$-algebra.
\end{proposition}
\demo
Let $f,g,l \in \Fr$ and $(x,h_t(x))\in R$. Then:
\begin{align*}
   & (f \ast g)\ast l(x,h_t(x)) \\
   & = \sum_{s \in G} (f \ast g)(x,h_s(x)) l(h_s(x),h_t(x)) \\
   & = \sum_{s \in G}\left( \sum_{r \in G}f(x,h_r(x))g(h_r(x),h_s(x)) \right) l(h_s(x), h_t(x))  \\
   & = \sum_{s \in G} \sum_{r \in G}[f(x,h_r(x))g(h_r(x),h_s(x))]l(h_s(x), h_t(x))\\
   & = \sum_{s \in G} \sum_{r\in G}f(x,h_r(x))[g(h_r(x),h_s(x))l(h_s(x), h_t(x))], \,\, \mbox{since $\K$ is associative} \\
   & = \sum_{r \in G} f(x,h_r(x)) \sum_{s \in G} g(h_r(x),h_s(x)) l(h_s(x),h_t(x))\\
   & = \sum_{r \in G} f(x,h_r(x)) \sum_{s \in G} g \ast l (h_r(x),h_t(x)) \\
  & =  f\ast ( g\ast l) (x,h_t(x)).
\end{align*}
So $(f\ast g)\ast l = f\ast ( g\ast l)$ and $\Fr$ is associative.
\cqd

\begin{corolario}
$\Fo \rtimes_{\alpha}G$  is an associative $\K$-algebra.
\end{corolario}

To finalize the paper we use theorem \ref{oteorema} to characterize the ideals in $\Fo \rtimes_{\alpha}G$ in terms of $R$-invariant sets:

\begin{definition}
A subset $Z\subseteq X$ is said $R$-invariant if whenever $(z,x)\in R$, with $z\in Z$, then $x\in Z$.
\end{definition}

\begin{proposition}
Let $I$ be an ideal in $\Fr$. Then there exists a $R$-invariant set $Z$ such that $I= \mathcal{F}_0\left( (Z\times Z)\cap R \right)$.
\end{proposition}
\demo
First notice that if $f$ is a non zero function in $I$ and $f(y,h_t(y))\neq 0$ then $\delta_{(y,h_t(y))} \in I$. This follows because $\delta_{(y,h_t(y))} = \delta_{(y,y)}\ast f \ast \delta_{(h_t(y),h_t(y))}$. Also notice that $\delta_{(y,h_t(y))} \ast \delta_{(z,h_s(z))} = 0$ if $z\neq h_t(y)$ and $\delta_{(y,h_t(y))} \ast \delta_{(z,h_s(z))} = \delta_{(y,h_{st}(y))}$ if $z= h_t(y)$. 

With the above in mind, let $I$ be an ideal in $\Fr$. Define $$Z=\{z\in X: \exists \ f \in I \text{ and }  t\in G \ \text{ such that } f(z,h_t(z)) \neq 0\}.$$

We will prove first that $Z$ is $R$-invariant. For this, let $z \in Z$ and $(z, h_s(z))\in R$. Then, by the definition of $Z$, there exists $f\in I$ and $t\in G$ such that $f(z,h_t(z))\neq 0$ and, by the first paragraph in this proof, we have that $\delta_{(z,h_t(z))} \in I$. But then $\delta_{(z,h_t(z))} \ast \delta_{(h_t(z),h_{st^{-1}}(z))}=\delta_{(z,h_s(z))} \in I$. We conclude that $\delta_{(h_s(z),z)} \ast \delta_{(z,h_s(z))} = \delta_{(h_s(z),h_s(z))} \in I$ and hence $h_s(z) \in Z$ as desired.

Next we prove that $\mathcal{F}_0\left( (Z\times Z)\cap R \right) \subseteq I$. Notice that for this it is enough to show that for each $(x,h_t(x)) \in (Z\times Z)\cap R$ the associated delta Dirac function $\delta_{(x,h_t(x))}$ belongs to $I$. So, let $(x,h_t(x)) \in (Z\times Z)\cap R$. Then, since $x\in Z$, there exists $f\in I$ and $s\in G$ such that $f(x,h_s(x))\neq 0$ and hence $\delta_{(x,h_s(x))}\in I$. We conclude that $\delta_{(x,h_t(x))} = \delta_{(x,h_s(x))}\ast \delta_{(h_s(x),h_{ts^{-1}}(x))}$ belongs to $I$ as desired.

Finally, $I\subseteq \mathcal{F}_0\left( (Z\times Z)\cap R \right)$, since if $f\in I$ and $f((x,h_t(x)))\neq 0 $ then $z\in Z$ and hence, by the $R$-invariance of $Z$, $h_t(z)\in Z$. So the support of $f$ is contained in $Z\times Z$.
\cqd 

\begin{obs} Notice that if $Z$ is an $R$-invariant subset then $\mathcal{F}_0\left( (Z\times Z)\cap R \right)$ is an ideal of $\Fr$.
\end{obs}

\begin{obs} Recently Leavitt path algebras of countable graphs were characterized as partial skew group rings, see \cite{GR1}. In the case of finite graphs with no cycles the partial skew group ring introduced in \cite{GR1} arises from a free partial action on a set. So, for these graphs, one can apply the last result above to obtain the characterization of the ideals of the associated Leavitt path algebras given by Tomforde in \cite{tomforde}. We refrain to provide the details as this would require the introduction of many notions not mentioned here.
\end{obs}

\vspace{1.5pc}

V. M. Beuter, Departamento de Matemática, Universidade do Estado de Santa Catarina, Joinville, 89219-710, Brasil

Email: vivibeuter@gmail.com

\vspace{0.5pc}

D. Goncalves, Departamento de Matemática, Universidade Federal de Santa Catarina, Florianópolis, 88040-900, Brasil

Email: daemig@gmail.com

\vspace{0.5pc}

\end{document}